\input amstex
\documentstyle{amsppt}
%
%
\nopagenumbers

\def\grad{\operatorname{grad}}
\def\Ker{\operatorname{Ker}}
\def\Cl{\operatorname{Cl}}
\accentedsymbol\tih{\tilde h}
\accentedsymbol\tiW{\tilde W}
\def\negskp{\hskip -2pt}
\def\compos{\,\raise 1pt\hbox{$\sssize\circ$} \,}
\pagewidth{360pt}
\pageheight{606pt}
\leftheadtext{Ruslan A. Sharipov}
\rightheadtext{First problem of globalization \dots}
\topmatter
\title First problem of globalization\\
{\lowercase{in the theory of dynamical systems
\vadjust{\vskip -1.0ex}admitting\\ the normal
shift of hypersurfaces.}}
\endtitle
\author
R.~A.~Sharipov
\endauthor
\abstract
Formula for the force field of Newtonian dynamical systems
admitting the normal shift of hypersurfaces in Riemannian
manifolds is considered. Problem of globalization for
geometric structures associated with this formula is studied.
\endabstract
\address Rabochaya street 5, 450003, Ufa, Russia
\endaddress
\email \vtop to 20pt{\hsize=280pt\noindent
R\_\hskip 1pt Sharipov\@ic.bashedu.ru\newline
ruslan-sharipov\@usa.net\vss}
\endemail
\urladdr
http:/\negskp/www.geocities.com/CapeCanaveral/Lab/5341
\endurladdr
\endtopmatter
\loadbold
\document
\head
1. Introduction.
\endhead
     Let $M$ be a Riemannian manifold of the dimension $n$. Newtonian
dynamical system in $M$ is determined in local coordinates by $n$
ordinary differential equations
$$
\ddot x^k+\sum^n_{i=1}\sum^n_{j=1}\Gamma^k_{ij}\,\dot x^i\,\dot x^j
=F^k(x^1,\ldots,x^n,\dot x^1,\ldots,\dot x^n),\hskip -2em
\tag1.1
$$
where $k=1,\,\ldots,\,n$. Here $\Gamma^k_{ij}=\Gamma^k_{ij}(x^1,\ldots,
x^n)$ are components of metric connection, while $F^k$ are components
of force vector $\bold F$. They determine force field of dynamical
system \thetag{1.1}. Let $S$ be a hypersurface in $M$ and let $p\in M$. 
Consider the following initial data for the system of equations
\thetag{1.1}:
$$
\xalignat 2
&\quad x^k\,\hbox{\vrule height 8pt depth 8pt width 0.5pt}_{\,t=0}
=x^k(p),
&&\dot x^k\,\hbox{\vrule height 8pt depth 8pt width 0.5pt}_{\,t=0}
=\nu(p)\cdot n^k(p).\hskip -2em
\tag1.2
\endxalignat
$$
Here $n^k(p)$ are components of unitary normal vector $\bold n$ to
$S$ at the point $p$. Initial data \thetag{1.2} determine the trajectory
coming out from the point $p$ in the direction of normal vector $\bold
n(p)$. The quantity $\nu(p)$ in \thetag{1.2} is introduced to determine
modulus of initial velocity for such trajectory.\par
    Let's choose and fix some point $p_0\in S$, then consider a smooth
function $\nu(p)$ defined in some neighborhood of the point $p_0$. Let
$$
\nu(p_0)=\nu_0\neq 0.\hskip -2em
\tag1.3
$$
Then in some (possibly smaller) neighborhood of $p_0$ this function
$\nu(p)$ does not vanish and take values of some definite sign. Upon
restricting $\nu(p)$ to such neighborhood we use it to determine
initial velocity in \thetag{1.2}. As a result we obtain a family of
trajectories of dynamical system \thetag{1.1}. Displacement of points 
of hypersurface $S$ along these trajectories determines shift maps
$f_t\!:S\to S_t$. Relying upon theorem on smooth dependence on initial
data for the system of ODE's (see \cite{1} and \cite{2}), we can assume
that shift maps $f_t\!:S'\to S'_t$ are defined in some neighborhood
$S'$ of the point $p_0$ on $S$ for all values of parameter $t$ in some
interval $(-\varepsilon,\,+\varepsilon)$ on real axis $\Bbb R$.
At the expense of further restriction of the interval $(-\varepsilon,\,
+\varepsilon)$ one can make maps $f_t\!:S'\to S'_t$ to be diffeomorphisms
and make their images $S'_t$ to be smooth hypersurfaces, disjoint union
of which fills some neighborhood of the point $p_0$ in $M$. Moreover,
at the expense of restriction of the neighborhood $S'$ and the range
of parameter $t$ one can reach the situation in which shift trajectories
would cross hypersurfaces $S_t$ transversally at all points of mutual
intersection. For such a case we state the following definitions.
\definition{Definition 1.1} Shift $f_t\!:S'\to S'_t$ of some part $S'$
of hypersurface $S$ along trajectories of Newtonian dynamical system
\thetag{1.1} is called {\bf a normal shift} if all hypersurfaces $S'_t$
arising in the process of shifting are perpendicular to the trajectories
of this shift.
\enddefinition
\definition{Definition 1.2} Newtonian dynamical system \thetag{1.1}
with force field $\bold F$ ia called {\bf a system admitting normal
shift} in strong sense if for any hypersurface $S$ in $M$, for any
point $p_0\in S$, and for any real number $\nu_0\neq 0$ one can find
a neighborhood $S'$ of the point $p_0$ on $S$, and a smooth function
$\nu(p)$, which do not vanish in $S'$ and which is normalized by the
condition \thetag{1.3}, such that the shift $f_t\!:S'\to S'_t$ defined
by this function is a normal shift in the sense of definition~1.1. 
\enddefinition
    First we used the definition without normalizing condition
\thetag{1.3} for the function $\nu(p)$. Such definition was called
the normality condition. Definition~1.2 strengthens this condition
making it more restrictive with respect to the choice of force
field $\bold F$ of dynamical system \thetag{1.1}. Therefore it is
called strong normality condition.\par
     Definitions~1.1 and 1.2 form the base of the theory of dynamical
systems admitting the normal shift. This theory was constructed in
papers \cite{3--18}. The results of these papers were used in preparing
theses \cite{19} and \cite{20}.\par
     As it was shown in \cite{19}, Newtonian dynamical systems admitting
the normal shift of hypersurfaces in Riemannian manifolds of the dimension
$n\geqslant 3$ can be effectively described. Force field of such systems
is given by explicit formula
$$
F_k=\frac{h(W)\,N_k}{W_v}-v\sum^n_{i=1}\frac{\nabla_iW}
{W_v}\,\bigl(2\,N^i\,N_k-\delta^i_k\bigr),\hskip -2em
\tag1.4
$$
which contains one arbitrary function of one variable $h=h(w)$ and
one arbitrary function of $(n+1)$ variables $W=W(x^1,\ldots,x^n,v)$
restricted by natural condition
$$
W_v=\frac{\partial W}{\partial v}\neq 0.\hskip -2em
\tag1.5
$$
Components of gradient $\nabla W$ in formula \thetag{1.4} are partial
derivatives
$$
\nabla_iW=\frac{\partial W}{\partial x^i}.\hskip -2em
\tag1.6
$$
Here $N^i$ and $N_k$ are components of unitary vector
$\bold N$ directed along velocity vector:
$$
\pagebreak
\xalignat 2
&N^i=\frac{v^i}{|\bold v|},&&N_k=\frac{v_k}{|\bold v|}.\hskip -2em
\tag1.7
\endxalignat
$$
Upon substituting \thetag{1.5}, \thetag{1.6}, and \thetag{1.7} into
formula \thetag{1.4} independent variable $v$ should be replaced by
modulus of velocity vector: $v=|\bold v|$.
\head
2. First problem of globalization.
\endhead
    If we fix a pair of functions $(h,W)$, then formula \thetag{1.4}
uniquely determines the force field $\bold F$ of Newtonian dynamical
system \thetag{1.1}. However, fixing force field \thetag{1.4}, we
cannot determine uniquely the corresponding pair of functions $(h,W)$.
In particular, global force field $\bold F$ can be represented by
different pairs of functions in different local maps forming an atlas
of the manifold $M$. Thus, we meet a problem of describing global
geometric structures associated with such a way of defining force
field $\bold F$. This problem was formulated by S.~E.~Kozlov and
Yu\.~R.~Romanovsky when I was reporting my thesis \cite{19} in the
seminar of N.~Yu\.~Netsvetaev at Saint-Petersburg department of
Steklov mathematical Institute.\par
    There is another problem of globalization concerning the process
of normal shift of some particular hypersurface $S$ along trajectories
of dynamical system \thetag{1.4}. We shall call it second problem of
globalization, though historically it arises earlier than first one.
Second problem was formulated by A.~S.~Mishchenko when I was reporting
thesis \cite{19} in the seminar of the Chair of higher geometry and
topology at Moscow State University. It's expedient to deal with
second problem of globalization only upon solving first one. Therefore
we shall consider it in separate paper.
\head
3. Scalar ansatz and gauge transformations.
\endhead
    Let's consider the projection of force vector \thetag{1.4} onto
the direction of the velocity vector. This projection can be calculated
as a scalar product of vectors $\bold F$ and $\bold N$:
$$
A=(\bold F\,|\,\bold N)=\sum^n_{k=1}F_k\,N^k.\hskip -2em
\tag3.1
$$
Substituting \thetag{1.4} into \thetag{3.1}, we get the following
expression for $A$:
$$
A=\frac{h(W)}{W_v}-\frac{v}{W_v}\,(\nabla W\,|\,\bold N).\hskip -2em
\tag3.2
$$
A very important point is that force fields \thetag{1.4} can be recovered
by corresponding scalar fields $A$. This recovery is given by formula,
which is called {\bf a scalar ansatz}:
$$
F_k=A\,N_k-|\bold v|\,\sum^n_{i=1}P^i_k\,\tilde\nabla_i A.
\hskip -2em
\tag3.3
$$
Here $P^i_k=\delta^i_k-N^i\,N_k$ are components of orthogonal projector
onto the hyperplane perpendicular to the vector $\bold v$. By $\tilde
\nabla_i A$ we denote partial derivatives $\partial A/\partial v^i$.
Scalar ansatz \thetag{3.3} was found in \cite{18}. In thesis \cite{19}
it was used in deriving formula \thetag{1.4}.\par
     Formulas \thetag{3.1} and \thetag{3.3} set up one-to-one correspondence
of vector fields $\bold F$ of the form \thetag{1.4} with scalar fields
$A$ of the form \thetag{3.2}. Formula \thetag{3.2} uniquely determines
the scalar field $A$ by the pair of functions $(h,W)$. But inverse
correspondence is not univalent. This is confirmed by the existence of
gauge transformations
$$
\aligned
&W(x^1,\ldots,x^n,v)\longrightarrow \rho(W(x^1,\ldots,x^n,v)),\\
&h(w)\longrightarrow h(\rho^{-1}(w))\cdot\rho'(\rho^{-1}(w))
\endaligned\hskip -2em
\tag3.4
$$
with one arbitrary function of one variable $\rho=\rho(w)$.
Transformations \thetag{3.4} change $h$ and $W$, but they don't change
the scalar field $A$.\par
     Let's investigate which part of information on $h$ and $W$ can be
recovered by $A$. Suppose that the point $p\in M$ is fixed. The
dependence of $A$ on the direction of velocity vector at the point $p$
is determined by the term $\bold N$ in scalar product
$(\nabla W\,|\,\bold N)$. Therefore if we change $\bold v$ by $-\bold v$,
first summand in \thetag{3.2} remains unchanged, while second summand
changes in sign. Hence
$$
\xalignat 2
&\frac{h(W)}{W_v}=\frac{A(\bold v)+A(-\bold v)}{2},
&&\frac{(\nabla W\,|\,\bold N)}{W_v}=\frac{A(-\bold v)
-A(\bold v)}{2\,|\bold v|}.
\endxalignat
$$
Keeping the value of $v=|\bold v|$ unchanged, we can change the direction
of vector $\bold N$. This allows us to determine each component of vector
$\nabla W/W_v$. Thus by $A$ one can recover the scalar $h(W)/W_v$ and the
vector $\nabla W/W_v$.\par
    Let $p$ be a point of the manifold $M$. Suppose that the field $A$
is determined by two pairs of functions $(h,W)$ and $(\tih,\tiW)$ in some
neighborhood of $p$. Then
$$
\xalignat 2
&\quad\frac{h(W)}{W_v}=\frac{\tih(\tiW)}{\tiW_v},
&&\frac{\nabla W}{W_v}=\frac{\nabla\tiW}{\tiW_v}.
\hskip -2em
\tag3.5
\endxalignat
$$
Being more accurate, one should note that functions $W$ and $\tiW$ are
determined in some domain $U$ in Cartesian product $M\times\Bbb R^{\sssize
+}$, where by $\Bbb R^{\sssize +}$ we denote the set of positive numbers.
Second relationship in \thetag{3.5} means that complete gradients of these
two functions in $U$ are collinear:
$$
\grad W\parallel\grad\tiW.\hskip -2em
\tag3.6
$$
The conditions $W_v\neq 0$ and $\tiW_v\neq 0$ mean that both gradients
in \thetag{3.6} are nonzero. This situation is described by the following
lemma.
\proclaim{Lemma 3.1} If gradient of one smooth function $f(x^1,\ldots,x^n)$
is nonzero in some domain $U\subset\Bbb R^n$ and gradient of another
smooth function $g(x^1,\ldots,x^n)$ is collinear to it in $U$, then
functions $f$ and $g$ are functionally dependent in $U$. This means that
for each point $p\in U$ one can find some neighborhood $O(p)$ and a smooth
function of one variable $\rho(y)$ such that $g=\rho\,\compos f$ in $O(p)$.
\endproclaim
    Lemma~3.1 is purely local fact following from the theory of implicit
functions (see \cite{21} and \cite{22}). But, despite to this, it is worth
while, since it describes the structure of non-uniqueness in inverse
correspondence for $(h,W)\to A$.
\proclaim{Theorem 3.1} Suppose that two pairs of functions $(h,W)$ and
$(\tih,\tiW)$ defined in some domain $U\subset M\times\Bbb R^{\sssize +}$
determine the same force field $\bold F$ of the form \thetag{1.4}. Then
for each point $q\in U$ one can find some neighborhood $O(q)$ and a smooth
function of one variable $\rho(y)$ such that $(h,W)$ and $(\tih,\tiW)$ are
bound by the gauge transformation \thetag{3.4} in $O(q)$.
\endproclaim
\head
4. Projectivization of cotangent bundle.
\endhead
     Denote by $\Cal M$ the Cartesian product of manifolds
$M\times\Bbb R^{\sssize +}$. Let $\Cal T^*\!\Cal M$ be cotangent
bundle for $\Cal M$. If we take the pair of functions $h$ and $W$,
\pagebreak which determine force field $\bold F$ of the form
\thetag{1.4}, then we see that derivatives
$$
\nabla_1W,\ \ \nabla_2,\ \ \ldots,\ \ \nabla_nW,\ \ W_w
$$
constitute the set of components of differential 1-form $\boldsymbol
\omega=dW$. The domain, where this 1-form is defined, shouldn't
coincide with the whole manifold $\Cal M$. Hence we have a local
section of the bundle $\Cal T^*\!\Cal M$. Second summand in formula
\thetag{1.4} contain not the components of differential form
$\boldsymbol\omega$ by themselves, but the quotients
$$
b_i=-\frac{\nabla_iW}{W_v}=-\frac{\omega_i}{\omega_{n+1}}.\hskip -2em
\tag4.1
$$
Let's factorize fibers of cotangent bundle $\Cal T^*\!\Cal M$ with
respect to the action of multiplicative group of real numbers
$\boldsymbol\omega\to\alpha\cdot\boldsymbol\omega$. In other words,
we replace linear spaces $\Cal T^*_q(\Cal M)$ over the points
$q\in\Cal M$ by corresponding projective spaces $\Cal P^*_{\!q}(\Cal M)$. 
As a result we get projectivized cotangent bundle $\Cal P^*\!\Cal M$.
This is locally trivial bundle $\Cal P^*\!\Cal M$, standard fiber of
which is $n$-dimensional projective space $\Bbb R\Bbb P^n$ (see
definitions in books \cite{23} or \cite{24}).\par
     Fibers of projective bundle $\Cal P^*\!\Cal M$ are parameterized
by components of covectors $\boldsymbol\omega$ taken up to an arbitrary
numeric factor:
$$
\alpha\cdot\omega_1,\ \ \alpha\cdot\omega_2,\ \ \ldots,\ \
\alpha\cdot\omega_n,\ \ \alpha\cdot\omega_{n+1}.\hskip -2em
\tag4.2
$$
If $\omega_{n+1}\neq 0$, then we can choose numeric factor $\alpha=
1/\omega_{n+1}$. Then from \thetag{4.2} we obtain $-b_1,\ -b_2,\ \ldots,
\ -b_n,\ 1$. This means that quantities $b_i$ from \thetag{4.1} are the
local coordinates in one of affine maps in projective fiber of the
bundle $\Cal P^*\!\Cal M$. Let's turn back to the problem of globalization
formulated in section~2. From formulas \thetag{3.5} we derive the following
proposition.
\proclaim{Lemma 4.1} Each force field $\bold F$ of the form
\thetag{1.4} determines some global section $\sigma$ of projectivized
cotangent bundle $\Cal P^*\!\Cal M$.
\endproclaim
    But not all global sections of the bundle $\Cal P^*\!\Cal M$ can
be obtained in this way. There is a restriction. The matter is that
on the level of cotangent bundle $\Cal T^*\!\Cal M$ our section $\sigma$
in lemma~4.1 is represented by closed differential forms $\boldsymbol
\omega$, which possibly may be defined only locally. Let's study how
this fact is reflected on the level of projective bundle $\Cal P^*\!
\Cal M$\,? In order to recover components of the form $\boldsymbol\omega$
in \thetag{4.2} by $b_1,\ b_2,\ \ldots,\ b_n$ we should take a proper
factor $\varphi=\omega_{n+1}$. Then
$$
\omega_i=\cases -b_i\,\varphi&\text{for \ }i=1,\,\ldots,\,n,\\
\quad\varphi&\text{for \ }i=n+1.\endcases\hskip -2em
\tag4.3
$$
Closedness of the form $\boldsymbol\omega$ is written in form of
the following relationships:
$$
\frac{\partial\omega_i}{\partial x^j}-\frac{\partial\omega_j}
{\partial x^i}=0.\hskip -2em
\tag4.4
$$
Here we denote $v=x^{n+1}$. This is natural, since $\Cal M=M\times
\Bbb R^{\sssize +}$. Substituting \thetag{4.3} into the relationships
\thetag{4.4}, for $i\leqslant n$ and $j\leqslant n$ we get
$$
\frac{\partial b_i}{\partial x^j}\,\varphi+
\frac{\partial\varphi}{\partial x^j}\,b_i=
\frac{\partial b_j}{\partial x^i}\,\varphi+
\frac{\partial\varphi}{\partial x^i}\,b_j.
\hskip -2em
\tag4.5
$$
From the same relationships \thetag{4.4} for the case $i\leqslant n$ and
$j=n+1$ we derive
$$
\frac{\partial\varphi}{\partial x^i}=-\frac{\partial b_i}{\partial v}\,
\varphi-\frac{\partial\varphi}{\partial v}\,b_i.\hskip -2em
\tag4.6
$$
Now let's substitute the derivatives $\partial\varphi/\partial x^i$ and
$\partial\varphi/\partial x^j$ calculated according to \thetag{4.6} into
the equations \thetag{4.5}. As a result we obtain the equations free of
$\varphi$:
$$
\left(\frac{\partial}{\partial x^j}+b_j\,\frac{\partial}
{\partial v}\right)b_i=\left(\frac{\partial}{\partial x^i}
+ b_i\,\frac{\partial}{\partial v}\right)b_j.\hskip -2em
\tag4.7
$$
Note that the equations \thetag{4.7} are already known (see \cite{19},
Chapter~\uppercase\expandafter{\romannumeral 7}, \S\,4). However, the
geometric interpretation of quantities $b_i$ in \cite{19} was quite
different.\par
\proclaim{Lemma 4.2} Each force field $\bold F$ of the form \thetag{1.4}
determines some global section $\sigma$ of the bundle $\Cal P^*\!\Cal M$
with components satisfying the equations \thetag{4.7}.
\endproclaim
     The equations \thetag{4.7} above arise as necessary condition for
the existence of closed differential 1-form $\boldsymbol\omega$
corresponding to the section of projective bundle $\Cal P^*\!\Cal M$.
But they are sufficient condition for the existence of such 1-form
as well (certainly, only for local existence). Let's prove this fact.
In order to integrate the equations \thetag{4.6} we use the auxiliary
system of Pfaff equations
$$
\frac{\partial V}{\partial x^i}=b_i(x^1,\ldots,x^n,V)
\text{, \ where \ }i=1,\,\ldots,\,n.\hskip -2em
\tag4.8
$$
The relationships \thetag{4.7} are exactly the compatibility conditions
for the equations \thetag{4.8}. Remember that variables $x^1,\,\ldots,
\,x^n,\,v$ are local coordinates in the manifold $\Cal M=M\times\Bbb
R^{\sssize +}$, while first $n$ of them are local coordinates in
$M$. Let's fix some point $p_0\in M$. Without loss of generality we
can assume that local coordinates of the point $p_0$ are equal to zero.
For compatible system of Pfaff equations \thetag{4.8} we set up the
following Cauchy problem at the point $p_0$:
$$
V\,\hbox{\vrule height 8pt depth 8pt width 0.5pt}_{\,x^1=\,\ldots\,
\,=\,x^n\,=\,0}=w.\hskip -2em
\tag4.9
$$
Thereby we take $w>0$. The solution of Cauchy problem \thetag{4.9}
for the equations \thetag{4.8} does exist and it is unique in some
neighborhood of the point $p_0$. It is smooth function of coordinates
$x^1,\,\ldots,\,x^n$ and parameter $w$:
$$
v=V(x^1,\ldots,x^n,w).\hskip -3em
\tag4.10
$$
For $x^1=\ldots=x^n=0$ due to \thetag{4.9} we have $V(0,\ldots,0,w)=w$.
Therefore
$$
\frac{\partial V}{\partial w}\,\hbox{\vrule height 14pt depth
10pt width 0.5pt}_{\,x^1=\,\ldots\,\,=\,x^n\,=\,0}=1.\hskip -3em
\tag4.11
$$
Let's consider the set of point $q=(p_0,v)$ in $\Cal M$. They form a
linear ruling in Cartesian product $\Cal M=M\times\Bbb R^{\sssize +}$.
Let's denote it by $l_0=l(p_0)$. \pagebreak The equality \thetag{4.11} means
that for any point $q_0\in l_0$ there is some neighborhood of this
point, where we have local coordinates $y^1,\,\ldots,\,y^n,\,w$ related
to initial coordinates $x^1,\,\ldots,\,x^n,\,v$ as 
$$
\cases x^i=y^i\text{\ \ for \ }i=1,\,\ldots,\,n,\\
v=V(y^1,\ldots,y^n,w).\endcases\hskip -3em
\tag4.12
$$
Back transfer to initial coordinates is determined by a function
$W(x^1,\ldots,x^n,v)$:
$$
\cases y^i=x^i\text{\ \ for \ }i=1,\,\ldots,\,n,\\
w=W(x^1,\ldots,x^n,v).\endcases\hskip -3em
\tag4.13
$$
Function $W(x^1,\ldots,x^n,v)$ is calculated implicitly from the
relationship \thetag{4.10} considered as the equation with respect
to $w$.\par
    Let's use \thetag{4.12} and \thetag{4.13} for to simplify the
equations \thetag{4.6}. Instead of function $\varphi(x^1,\ldots,x^n,v)$
in these equations we introduce another function
$$
\psi(y^1,\ldots,y^n,w)=\varphi(y^1,\ldots,y^n,V(y^1,\ldots,y^n,w)).
\hskip -3em
\tag4.14
$$
The equations \thetag{4.6} are reduced to the following equations
for the function \thetag{4.14}:
$$
\frac{\partial\psi}{\partial y^i}=-B_i\,\psi.\hskip -3em
\tag4.15
$$
The quantities $B_i$ are expressed through the derivatives of the
function $V$:
$$
B_i=\frac{1}{Z}\,\frac{\partial Z}{\partial y^i}\text{, \ where \ }
Z=\frac{\partial V}{\partial w}.\hskip -3em
\tag4.16
$$
It is easy to see that \thetag{4.15} is a system of Pfaff equations
being compatible due to \thetag{4.16}. Moreover, it is explicitly
integrable. General solution of the equations \thetag{4.15} is given
by the following explicit formula:
$$
\psi=\frac{C(w)}{Z(y^1,\ldots,y^n,w)},\hskip -3em
\tag4.17
$$
Here $C(w)$ is an arbitrary smooth function of one variable. Now let's
use the local invertibility of the relationship \thetag{4.14}:
$$
\varphi(x^1,\ldots,x^n)=\psi(x^1,\ldots,x^n,W(x^1,\ldots,x^n,v)).
\hskip -3em
\tag4.18
$$
From \thetag{4.17} and \thetag{4.18} we derive general solution for
the system of equations \thetag{4.6}:
$$
\varphi=C(W)\cdot W_v\text{, \ where \ }W_v=\frac{\partial W}{\partial v}.
$$
Similar to force field $\bold F$ in formula \thetag{1.4}, it is
determined by two functions $C(w)$ and $W(x^1,\ldots,x^n,v)$, the
latter one satisfying the condition \thetag{1.5}. This coincidence
is not occasional. From \thetag{4.8} and from \thetag{4.18} for
$b_i$ we derive the relationship
$$
\pagebreak 
b_i=-\frac{\nabla_i W}{W_v},\hskip -3em
\tag4.19
$$
being of the same form as \thetag{4.1}. Certainly, function $W$
in \thetag{4.19} obtained by inverting local change of variables
\thetag{4.12} shouldn't coincide with initial function $W$ in
\thetag{4.1}. The relation of these two functions is characterized
by theorem~3.1 (see above). The calculations we have just made result
in the following lemma, sharpening lemma~4.2.
\proclaim{Lemma 4.3} The relationships \thetag{4.7} form necessary
and sufficient condition for global section $\sigma$ of projectivized
cotangent bundle $\Cal P^*\!\Cal M$ given by its components
$b_1,\,\ldots,\,b_n$ in local coordinates to be related to some
force field $\bold F$ of the form \thetag{1.4}.
\endproclaim
\head
5. Involutive distributions.
\endhead
    Relying upon lemmas~4.2 and 4.3, now we consider some global section
$\sigma$ of projectivized cotangent bundle $\Cal P^*\!\Cal M$ that
satisfies the equations \thetag{4.7}. Let's reveal invariant meaning
of these equations. In order to do it we consider vector fields 
$$
\bold L_i=\frac{\partial}{\partial x^i}+b_i\,\frac{\partial}{\partial v}
\text{, \ where \ }i=1,\,\ldots,\,n,\hskip -2em
\tag5.1
$$
and some differential 1-form $\boldsymbol\omega$ with components
\thetag{4.3}. Values of vector fields \thetag{5.1} are linearly
independent at each point of the domain, where they are defined.
These values belong to the kernel of the form $\boldsymbol\omega$
for any choice of function $\varphi$ in \thetag{4.3}. The equations
\thetag{4.7} are exactly the commutation conditions for vector
fields \thetag{5.1}:
$$
[\bold L_i,\,\bold L_j]=0.\hskip -2em
\tag5.2
$$
Note that global sections of the bundle $\Cal P^*\!\Cal M$ are in
one-to-one correspondence with $n$-dimensional distributions in
the manifold $\Cal M$, the dimension of which is equal to $n+1$.
Indeed, in the neighborhood of each point $q\in\Cal M$ the section
$\sigma$ of the bundle $\Cal P^*\!\Cal M$ is determined by some
1-form $\boldsymbol\omega$ fixed up to a scalar factor $\varphi$.
But the kernel $U=\Ker\boldsymbol\omega\subset\Cal T_q(\Cal M)$
does not depend on this factor. Therefore we have global
$n$-dimensional distribution $U=\Ker\sigma$. And conversely, if
$n$-dimensional distribution $U$ is given, then in the neighborhood
of each point $q\in\Cal M$ we have 1-form $\boldsymbol\omega$
such that $U=\Ker\boldsymbol\omega$. Form $\boldsymbol\omega$ defines
local section of the bundle $\Cal P^*\!\Cal M$ in the neighborhood
of the point $q$. The fact that form $\boldsymbol\omega$ is determined
by $U$ uniquely up to a scalar factor means that local sections of
the bundle $\Cal P^*\!\Cal M$ are glued into one global section
$\sigma$ of this bundle.\par
     The condition \thetag{5.2} means that the distribution $U=\Ker
\sigma$ is involutive (see \cite{24}). In this case in the
neighborhood of each point $q\in\Cal M$ the section $\sigma$ can be
represented by closed 1-form $\boldsymbol\omega$. Let's introduce
the following terminology.
\definition{Definition 5.1} The section $\sigma$ of projectivized
cotangent bundle $\Cal P^*\!\Cal M$ is called {\bf closed} if
corresponding distribution $U=\Ker\sigma$ in $\Cal M$ is involutive.
\enddefinition
    For the sections $\sigma$ related to force fields \thetag{1.4}
the manifold $\Cal M$ is a Cartesian product $M\times\Bbb R^{\sssize +}$.
In this case we have a restriction expressed by the condition \thetag{1.5}.
It can be written as $\omega_{n+1}\neq 0$. Therefore we have the following
lemma.
\proclaim{Lemma 5.1} Global section $\sigma$ of projectivized cotangent
bundle $\Cal P^*\!\Cal M$ with base manifold $\Cal M=M\times\Bbb
R^{\sssize +}$ satisfies the condition \thetag{1.5} if and only if
corresponding distribution $U=\Ker\sigma$ is transversal to linear
rulings \pagebreak of cylindric manifold $M\times\Bbb R^{\sssize +}$.
\endproclaim
\noindent For the sake of brevity we shall write the condition stated
in lemma~5.1 as
$$
\Ker\sigma=U\nparallel\Bbb R^{\sssize +}.\hskip -2em
\tag5.3
$$
Results of lemmas~4.2, 4.3, and 5.1 can be summarized in the
following theorem.
\proclaim{Theorem 5.1} Each force field $\bold F$ of the form \thetag{1.4}
determines some closed global section $\sigma$ of projectivized cotangent
bundle $\Cal P^*\!\Cal M$ over base manifold $\Cal M=M\times\Bbb
R^{\sssize +}$ such that it satisfies additional condition \thetag{5.3}.
And conversely, each such section of the bundle $\Cal P^*\!\Cal M$
corresponds to some force field $\bold F$ of the form \thetag{1.4}.
\endproclaim
\head
6. Normalizing vector fields.
\endhead
    Up to now we studied only the second summand in formula \thetag{3.2}.
And we have found that it gives rise to geometric structures mentioned in
theorem~5.1. Now let's consider first summand in \thetag{3.2}. Denote by
$a$ the following quotient:
$$
a=\frac{h(W)}{W_v}.\hskip -2em
\tag6.1
$$
Function $a=a(x^1,\ldots,x^n,v)$ in \thetag{6.1} is invariant with
respect to gauge transformations \thetag{3.4}. Due to the relationships
\thetag{3.5} it can be continued through the region of overlapping of
two maps, in which force field $\bold F$ is determined by two different
pairs of functions $(h,W)$ and $(\tih,\tiW)$. But, despite to this fact,
it would be wrong to interpret $a$ as scalar field on $\Cal M$. The
matter is that in local coordinates, for which formula \thetag{1.4}
holds, the variable $v$ plays exclusive role related with the expansion
of $\Cal M$ into Cartesian product $M\times\Bbb R^{\sssize +}$. Due to
this reason we derive differential equations for the function $a=a(x^1,
\ldots,x^n,v)$. Let's apply one of the differential operators \thetag{5.1}
to $a$. This yields
$$
\left(\frac{\partial}{\partial x^i}+b_i\,\frac{\partial}{\partial v}
\right)a=h'(W)\,\frac{\nabla_iW+b_i\,W_v}{W_v}-\frac{h(W)}{W_v}\,
\frac{\nabla_iW_v+b_i\,W_{vv}}{W_v}.
$$
If we take into account \thetag{4.1}, then this relationship can be
brought to 
$$
\left(\frac{\partial}{\partial x^i}+b_i\,\frac{\partial}{\partial v}
\right)a=\frac{\partial b_i}{\partial v}\,a.\hskip -2em
\tag6.2
$$
Note that the equations \thetag{6.2} are also already known (see
\cite{19}, Chapter~\uppercase\expandafter{\romannumeral 7}, \S\,4).
Formula \thetag{1.4} was derived as a result of integrating the
equations \thetag{4.7} and \thetag{6.2}. Following \cite{19}, we
append vector fields \thetag{5.1} by the following one
$$
\bold L_{n+1}=a\,\frac{\partial}{\partial v}.\hskip -2em
\tag6.3
$$
The equations \thetag{6.2} are equivalent to the following commutational
relationships:
$$
[\bold L_i,\bold L_{n+1}]=0\text{, \ where \ }i=1,\,\ldots,\,n.\hskip -2em
\tag6.4
$$\par
    Now let's give invariant (coordianteless) interpretation for the
relationships \thetag{6.4}. Vector fields \thetag{5.1} by themselves
have no invariant interpretation. But their linear span at each point
$q$ coincides with $n$-dimensional subspace $U_q\subset \Cal T_q(\Cal M)$
defined by distribution $U=\Ker\sigma$. Consider one-dimensional
factor-spaces
$$
\varOmega_q=\Cal T_q(\Cal M)/U_q.\hskip -2em
\tag6.5
$$
They are glued into one-dimensional vector bundle $\varOmega\Cal M$ over
the base manifold $\Cal M=M\times\Bbb R^{\sssize +}$. Let $x^1,\,\ldots,
\,x^n,\,v$ be local coordinates in $\Cal M$ not necessarily related to
the structure of Cartesian product $M\times\Bbb R^{\sssize +}$, but such
that vector $\partial/\partial v$ is transversal to $U_q$. Then vectors
\thetag{5.1} form the base in subspace $U_q$, while elements of factor-space
\thetag{6.5} are cosets of subspace $U_q$ represented by vectors
\thetag{6.3}:
$$
a=\Cl_U(a\cdot\partial/\partial v).
$$
Sections of one dimensional vector bundle $\varOmega\Cal M$ in such local
coordinates can be associated with functions $a(x^1,\,\ldots,x^n,v)$ or
with vector fields 
$$
\bold X=a(x^1,\,\ldots,x^n,v)\cdot\frac{\partial}{\partial v}.
$$
\definition{Definition 6.1} Vector field $\bold X$ is called {\bf
normalizing field} for smooth distribution $U$ if for any vector
field $\bold Y$ belonging to $U$ the commutator $[\bold X,\,\bold Y]$
is also in $U$.
\enddefinition
\proclaim{Theorem 6.1} Smooth distribution $U$ of codimension $1$ in
the manifold $\Cal M$ possesses nonzero normalizing vector field
transversal to $U$ in the neighborhood of the point $q\in\Cal M$
if and only if it is involutive in the neighborhood of this point.
\endproclaim
\demo{Proof} Without loss of generality we can take $\dim\Cal M=n+1$
and $\dim U=n$. Since $\bold X\neq 0$, we can choose local
coordinates $x^1,\,\ldots,\,x^n,\,v$ in $\Cal M$ such that $\bold X
=\partial/\partial v$. And since $\bold X\nparallel U$, the base in
$U$ can be formed by vector fields $\bold L_i$ of the form \thetag{5.1}.
Let's write the condition that $\bold X$ is normalizing vector field
for the distribution $U$. For this purpose we calculate the commutator
$[\bold X,\,\bold L_i]$:
$$
[\bold X,\,\bold L_i]=-\frac{\partial b_i}{\partial v}\cdot
\frac{\partial}{\partial v}=-\frac{\partial b_i}{\partial v}
\cdot\bold X.
$$
Recall that vector field $\bold X$ is transversal to $U$. Therefore
from $[\bold X,\,\bold L_i]\in U$ we get
$$
\frac{\partial b_i}{\partial v}=0.\hskip -2em
\tag6.6
$$
If we take into account \thetag{6.6}, we find that the equations
\thetag{4.7} turn to identities. But we know, that they are equivalent
to commutational relationships \thetag{5.2}. Hence the distribution
$U$ is involutive. Theorem~6.1 is proved.
\qed\enddemo
    Let $\bold X$ be normalizing vector field for involutive distribution
$U$ and let $\bold Y$ be in $U$. Then $\bold X+\bold Y$ is also normalizing
vector field for $U$. Thus we can define normalizing sections of the
bundle $\varOmega\Cal M$ obtained by factorization of tangent bundle
$\Cal T\Cal M$ with respect to $U$.
\definition{Definition 6.2} Section $s$ of factor-bundle $\varOmega\Cal M=
\Cal T\Cal M/U$ is called {\bf normalizing section} if in the neighborhood
of each point $q\in\Cal M$ it is represented by some normalizing vector
\pagebreak field for the distribution $U$.
\enddefinition
     Now we can formulate the main result of this paper, characterizing
global geometric structures associated with formula \thetag{1.4} for the
force field $\bold F$. It follows from all what was said above.
\proclaim{Theorem 6.2} Defining Newtonian dynamical system admitting
the normal shift in Riemannian manifold $M$ is equivalent to defining
closed global section $\sigma$ for projectivized cotangent bundle
$\Cal P^*\!\Cal M$ with base $\Cal M=M\times \Bbb R^{\sssize +}$,
satisfying the condition $\Ker\sigma\nparallel\Bbb R^{\sssize +}$,
and to defining normalizing global section $s$ for one-dimensional
factor-bundle $\varOmega\Cal M=\Cal T\Cal M/U$, where $U=\Ker\sigma$.
\endproclaim
\head
7. Integration of geometric structures.
\endhead
    Formulating theorem~6.2, we have made a step forward in understanding
global geometry associated with formula \thetag{1.4} for the force field
$\bold F$. But as far as the effectiveness of calculations in coordinates
is concerned, we came back to situation, in which scalar field $A$ is
expressed by formula
$$
A=a+\sum^n_{i=1}b_i\,v^i,\hskip -2em
\tag7.1
$$
where quantities $a$ and $b_1,\,\ldots,\,b_n$ should be found as solutions
of the equations \thetag{4.7} and \thetag{6.2}. Formula \thetag{3.2} was
more effective. Therefore we have a natural question: {\it can one integrate
the equations \thetag{4.7} and \thetag{6.2} globally and find the pair of
functions $(h,W)$ that would define scalar field $A$ by formula \thetag{3.2}
and force field $\bold F$ by formula \thetag{1.4} on the whole  manifold
$\Cal M$?}\par
    According to theorem~6.2, each force field $\bold F$ of Newtonian
dynamical system admitting the normal shift is related with some unique
closed global section $\sigma$ of the bundle $\Cal P^*\!\Cal M$. If such
section is generated by closed global section $\boldsymbol\omega$ of
cotangent bundle $\Cal T^*\!\Cal M$, then we can construct the function
$W=W(q)$ on $\Cal M$ by integrating 1-form $\boldsymbol\omega$ along
the curve binding the point $q$ with some fixed point $q_0$ on $\Cal M$:
$$
W(q)=\int^{\,q}_{q_0}\boldsymbol\omega.\hskip -2em
\tag7.2
$$
Formula \thetag{7.2} yields the function $W(q)$ that possibly can be
multivalued, since first homotopic group $\pi_1(\Cal M)$ of the manifold
$\Cal M$ can be non-trivial. This ambiguity is admissible. It can be
eliminated by passing to universal covering of $\Cal M$.\par
    Apart from $\boldsymbol\sigma$, each force field $\bold F$ of
Newtonian dynamical system admitting the normal shift determines some
section of factor-bundle $\varOmega\Cal M=\Cal T\Cal M/U$, where
$U=\Ker\sigma$. Let's use the structure of Cartesian product $M\times
\Bbb R^{\sssize +}$ of $\Cal M$. This yields the vector field $\bold E$
directed along linear rulings in $\Cal M$. If $x^1,\,\ldots,\,x^n$ are
local coordinates in $M$ and if $v$ is natural variable ranging in
positive semiaxis $\Bbb R^{\sssize +}$, then in local coordinates
$x^1,\,\ldots,\,x^n,\,v$ in $\Cal M$ this field is given by formula
$\bold E=\partial/\partial v$. According to theorem~6.2, we have
$U=\Ker\sigma\nparallel\Bbb R^{\sssize +}$, i\.~e\. $U\nparallel\bold E$.
Therefore the section $s$ of the bundle $\varOmega\Cal M$ can be
represented by the vector field 
$$
\bold X=a\cdot\bold E.\hskip -2em
\tag7.3
$$
This representation is unique, coefficient $a$ in it is a scalar field
(a function) on $\Cal M$. The condition that $s$ is normalizing section
with respect to $U$ in local coordinates $x^1,\,\ldots,\,x^n,\,v$ is
expressed by the equations \thetag{6.2} for the function $a$. It's easy
to check that if $a$ satisfies the equations \thetag{6.2}, then the
function $\varphi=1/a$ satisfies the equations \thetag{4.6}. Hence if
section $s$ is nonzero at all points $q\in\Cal M$, then we can use
$\varphi=1/a$ as proper integrating factor in formula \thetag{4.3}
determining components of closed 1-form $\boldsymbol\omega$. Contracting
this form with vector field \thetag{7.3}, we get 
$$
\boldsymbol\omega(\bold X)=C(\boldsymbol\omega\otimes\bold X)=1.
\hskip -2em
\tag7.4
$$
Section $\sigma$ of the bundle $\Cal P^*\!\Cal M$ determines 1-form
$\boldsymbol\omega$ up to a scalar factor, formula \thetag{4.3}
fixes this factor within the domain of local coordinates $x^1,\,\ldots,
\,x^n,\,v$, while the condition \thetag{7.4} shows that 1-forms defined
locally by this procedure is glued into one global closed 1-form
$\boldsymbol\omega$. Substituting its components into \thetag{7.1}, we get
$$
A=\frac{1}{\omega_{n+1}}-\sum^n_{i=1}\frac{\omega_i\,v^i}{\omega_{n+1}}.
\hskip -2em
\tag7.5
$$
Scalar field \thetag{7.5} corresponds to the force field $\bold F$ with
components
$$
F_k=\frac{N_k}{\omega_{n+1}}-v\sum^n_{i=1}\frac{\omega_i}
{\omega_{n+1}}\,\bigl(2\,N^i\,N_k-\delta^i_k\bigr),\hskip -2em
\tag7.6
$$
\proclaim{Theorem 7.1} If the section $s$ of factor-bundle $\varOmega
\Cal M=\Cal T\Cal M/U$ corresponding to the force field $\bold F$ of
Newtonian dynamical system admitting the normal shift is nonzero at
all points $q\in\Cal M=M\times\Bbb R^{\sssize +}$, then there is
a global closed 1-form $\boldsymbol\omega$ determining $\bold F$
according to the formula \thetag{7.6}.
\endproclaim
    In thesis \cite{19} it was noted that if the function $h(w)$
in formula \thetag{1.4} is nonzero, then at the expense of gauge
transformation \thetag{3.4} one can make it identically equal to
unity. There this fact was understood as purely local. Theorem~7.1
shows that it is valid in global situation too.
\head
8. Absence of topological obstructions.
\endhead
    It is well-known that some geometric structures cannot be realized
in manifolds with non-trivial topology. Thus, on the sphere $S^2$ there
are no smooth vector fields without special points, where they
vanish. For geometric structures from theorem~6.2 we have no such
obstructions. Indeed, on any manifold $M$ one has a smooth function
$w=w(p)$ which is not identically zero. Let $W(p,v)=w(p)+v$, where
$v\in\Bbb R^{\sssize +}$. It's obvious that the function $W(p,v)$ on
Cartesian product $M\times\Bbb R^{\sssize +}$ satisfies the condition
\thetag{1.5}. This function defines some global force field $\bold F$
of the form \thetag{1.4} and all geometric structures from theorem~6.2
as well.
\head
9. Acknowledgements.
\endhead
    I am grateful to A.~S.~Mishchenko for the invitation to visit
Moscow and for the opportunity to report the results of thesis \cite{19}
and succeeding papers \cite{25}, \cite{26}, and \cite{27} in his
seminar at Moscow State University. I am grateful to N.~Yu\.~Netsvetaev
for the invitation to visit Saint-Petersburg and for the opportunity
to report the same results in the seminar at Saint-Petersburg department
of Steklov Mathematical Institute. I am especially grateful to
Yu\.~R.~Romanovsky who showed me some famous historic places of
Saint-Petersburg and talked me about them during my very short
visit to this wonderful city.\par 
    I am grateful to all participants of both seminars mentioned above
and to my colleague E.~G.~Neufeld from Bashkir State University for
fruitful discussions, which stimulated preparing this paper.\par
     This work is supported by grant from Russian Fund for Basic Research
(project No\nolinebreak\.~00\nolinebreak-01-00068, coordinator 
Ya\.~T.~Sultanaev), and by grant from Academy of Sciences of the
Republic Bashkortostan (coordinator N.~M.~Asadullin). I am grateful
to these organizations for financial support.

\enddocument
\end